\documentclass[journal]{IEEEtran}

\usepackage{mdwlist}
\usepackage[tbtags]{amsmath}
\usepackage{multirow}
\usepackage{graphicx}
\usepackage{subfigure}
\usepackage{rotating}
\usepackage{epstopdf}
\usepackage{empheq}
\usepackage{color}
\usepackage{caption}
\def\Plus{\texttt{+}}
\def\Minus{\texttt{-}}  
\usepackage{todonotes}

\begin{document}
\title{Impact of Forecast Errors on Expansion Planning of Power Systems with a Renewables Target}

\author{S.~Pineda,~\IEEEmembership{Member,~IEEE,}
        J.~M.~Morales,~\IEEEmembership{Member,~IEEE,}
        and~T.~B.~Boomsma% <-this % stops a space
\thanks{S. Pineda and T. B. Boomsma thank the Danish Council for Strategic Research for support through “5s - Future Electricity Markets” project, no. 12-132636/DSF.}% <-this % stops a space
\thanks{J. M. Morales thanks the Danish Council for Strategic Research for support through ENSYMORA project, no. 10-093904/DSF. J. M. Morales is also sponsored by DONG Energy}}% <-this % stops a space

% The paper headers
\markboth{Submitted to IEEE Transactions on Power Systems}%
{Shell \MakeLowercase{\textit{et al.}}: Bare Demo of IEEEtran.cls for Journals}

\maketitle

\begin{abstract}
This paper analyzes the impact of production  forecast errors on the expansion planning of a power system and investigates the influence of market design to facilitate the integration of renewable generation. For this purpose, we propose a stochastic programming modeling framework to determine the expansion plan that minimizes system-wide investment and operating costs, while ensuring a given share of renewable generation in the electricity supply. Unlike existing ones, this framework includes both a day-ahead and a balancing market so as to capture the impact of both production forecasts and the associated prediction errors. Within this framework, we consider two paradigmatic market designs that essentially differ in whether the day-ahead generation schedule and the subsequent balancing re-dispatch are co-optimized or not. The main features and results of the model set-ups are discussed using an illustrative four-node example and a more realistic 24-node case study.
\end{abstract}


\IEEEpeerreviewmaketitle

\section{Notation}

\subsection{Indexes and sets}
\begin{basedescript}{\desclabelstyle{\pushlabel}\desclabelwidth{2em}}
    \item[$i$] Index of units/lines/loads/projects.
    \item[$n$] Index of buses.    
    \item[$s$] Index of scenarios in the day-ahead market.
    \item[$r$] Index of scenarios in the balancing market.
    \item[$\mathcal{F}$] Set of transmission lines.    
    \item[$\mathcal{G}$] Set of generating units.
    \item[$\mathcal{G}_n$] Set of generating units connected to bus $n$.
    \item[$\mathcal{L}$] Set of loads.
    \item[$\mathcal{L}_n$] Set of loads connected to bus $n$.
    \item[$\mathcal{P}$] Set of available expansion projects.
    \item[$\mathcal{W}$] Set of renewable generating units.
\end{basedescript}

\subsection{Parameters}

\begin{basedescript}{\desclabelstyle{\pushlabel}\desclabelwidth{2em}}
    \item[$a_{in}$] Line-bus indicator equal to 1/-1 if bus $n$ is the sending/receiving bus of line $i\in\mathcal{F}$, and 0 otherwise.
    \item[$b_{i}$] Susceptance of line $i\in\mathcal{F}$ (p.u.).
    \item[$c_i$] Price offer for energy sale by device $i\in\mathcal{G} \cup \mathcal{L}$ in the day-ahead market (\$/MWh).
    \item[$c^{\Plus}_i$] Price offer for energy sale by device $i\in\mathcal{G} \cup \mathcal{L}$ in the balancing market (\$/MWh).
    \item[$c^{\Minus}_i$] Price offer for energy purchase by device $i\in\mathcal{G} \cup \mathcal{L}$ in the balancing market (\$/MWh).
    \item[$l_i$] Length of transmission line $i$ (miles).
    \item[$q^0_i$] Fixed investment cost of project $i\in\mathcal{P}$ (\$/year).
    \item[$q_i$] Variable investment cost of project $i\in\mathcal{P}$ (\$/MWyear).
    \item[$T$] Duration of the planning horizon (h).
    \item[$x_i$] Peak consumption of load $i\in\mathcal{L}$ (MW).
    \item[$x^{\rm max}_i$] Maximum capacity of project $i\in\mathcal{P}$ (MW).
    \item[$x^{\Plus}_i$] Capacity factor for energy sale of device $i\in\mathcal{G} \cup \mathcal{L}$ in the balancing market (\%).
    \item[$x^{\Minus}_i$] Capacity factor for energy repurchase of device $i\in\mathcal{G} \cup \mathcal{L}$ in the balancing market (\%).
    \item[$\pi_{s}$] Probability of scenario $s$.
    \item[$\pi_{sr}$] Probability of scenario $r$ conditional on the realization of scenario $s$.
    \item[$\widehat{\rho}_{is}$] Forecasted capacity factor of device $i\in\mathcal{G} \cup \mathcal{L} \cup \mathcal{F}$ in the day-ahead market (p.u.).
    \item[$\widetilde{\rho}_{isr}$] Realized capacity factor of device $i\in\mathcal{G} \cup \mathcal{L} \cup \mathcal{F}$ in the balancing market (p.u.).  
    \item[$\eta$] Minimum share of renewable generation (\%).
\end{basedescript}

\subsection{Variables}
\begin{basedescript}{\desclabelstyle{\pushlabel}\desclabelwidth{2em}}
    \item[$\widehat{p}_{is}$] Power dispatch of device $i\in\mathcal{G}\cup\mathcal{L}$ or power flow through device $i\in\mathcal{F}$ in day-ahead scenario $s$ (MW).    
    \item[$p^{\Plus}_{isr}$] Power dispatch increase of device $i\in\mathcal{G}\cup\mathcal{L}$ in balancing scenario $r$ conditional on the realization of scenario $s$ (MW).
    \item[$p^{\Minus}_{isr}$] Power dispatch decrease of device $i\in\mathcal{G} \cup \mathcal{L}$ in balancing scenario $r$ conditional on the realization of scenario $s$ (MW).
    \item[$u_i$] Binary variable equal to 1 if project $i\in\mathcal{P}$ is built, and 0 otherwise.
    \item[$x_i$] Capacity of project $i\in\mathcal{P}$ (MW).
    \item[$\widehat{\delta}_{ns}$] Voltage angle at bus $n$ in day-ahead scenario $s$  (rad).
    \item[$\widetilde{\delta}_{nsr}$] Voltage angle at bus $n$ in balancing scenario $r$ conditional on the realization of scenario $s$ (rad).
\end{basedescript}

\section{Introduction}\label{Intro}

With an aim to reduce greenhouse gas emissions, ensure adequacy of energy resources and guarantee the security of energy supply, many governments employ policy targets for the exploitation of renewable energy sources. As an example, the European Union has agreed to cover 20\% of its total energy consumption with renewable production by 2020, whereas some EU countries aim for even higher shares \cite{EC2009}. Such ambitious targets will inevitably demand a paradigm shift in the design of future electricity markets and the configuration of future power systems.

In contrast to conventional power sources, renewable production, such as wind and solar power, is characterized by being highly variable, largely unpredictable and non-dispatchable. For the purpose of power system expansion planning, it is therefore necessary to account for both the variability of stochastic production throughout the planning horizon and the forecast errors caused by the limited predictability of these energy sources \cite{Bludszuweit2008}. 

Power system expansion planning under uncertainty has long been subject to extensive study in the literature \cite{Latorre2003}. As an example, \cite{Gorenstein1993} considers the expansion of generating capacity under demand, cost and technical uncertainty. Likewise, the authors in \cite{Torre2008} determine the optimal expansion plan for the transmission network while taking into account uncertainty in the future demand. More recent contributions \cite{Baringo2012,Orfanos2013} further investigate how the transmission expansion planning of a power system is affected by variability in wind power production. However, these models disregard forecast errors of the stochastic generating units and thus, fail to capture the impact of uncertainty in wind power production on expansion planning decisions.

In this paper, we address the expansion planning of a power system with a significant share of renewable generation, considering both the variability and limited predictability of the stochastic production. We propose a stochastic programming framework to determine the expansion plan (for stochastic and dispatchable production capacity as well as for the transmission network) that minimizes investment and system operating costs, while ensuring a given target of renewable generation in the electricity supply. Unlike in existing ones, the operation costs computed in our modeling framework include two terms: a day-ahead dispatch cost, which depends on the forecast of stochastic production and demand, and a balancing cost, which is a function of the forecast errors and the flexibility provided by the conventional generating units. To the best of our knowledge, our model is the first to explicitly model the impact of production forecast errors on the optimal expansion planning of generation and transmission capacity.

To properly incorporate the impact of forecast errors, we require two electricity trading floors: the day-ahead market, which is cleared 24 to 36 hours in advance of system operation and the balancing market, which copes with real-time energy imbalances \cite{Kirschen2004}. Moreover, we consider in this paper two paradigmatic market designs, which represent two extremes of day-ahead market clearing in view of potential forecast errors in renewable production:
\begin{itemize}
\item[-] The first is an ideal market-clearing procedure that efficiently handles forecast errors by jointly optimizing the operation in the day-ahead and balancing markets  \cite{Pritchard2010, Morales2013}. Under this market design, our expansion model becomes a three-stage stochastic program with capacity expansion, day-ahead generation scheduling and balancing re-dispatch being first-stage, second-stage and third-stage decisions, respectively \cite{birge2011}.
\item[-] The second represents an inefficient market design where the day-ahead generation schedule and the subsequent balancing re-dispatch are not jointly optimized. Under this market design, our expansion model is likewise a three-stage stochastic program. However, to capture the sequential and non-cooptimized clearings of the day-ahead and balancing markets, we formulate the model as a bi-level program \cite{dempe2002}.
\end{itemize}

We adopt the view of a central planner that minimizes the costs of power system expansion and operation. It is known, though, that the expansion plan promoted by a central planner is equivalent to that induced by an electricity market under perfect competition at the investment and operational stages \cite{Smeers2011}. For investment models that account for imperfect competition, the reader is referred to \cite{Murphy2005, Nanduri2009, Wogrin2013}. It should be remarked, however, that due to the complexity of these models, they often not allow for the inclusion of uncertainty.

We consider here a static approach to expansion planning where the optimal capacities of generation and transmission are determined for a single representative year. We thereby implicitly assume that the central planner aims to determine an optimal future power system configuration, rather than establishing when each expansion project should be carried out \cite{Latorre2003}. Nevertheless, the proposed models can be adapted to take into account not only sizing and placement of investments, but also timing considerations, albeit considerably increasing their computational burden.

To facilitate computations, we assume known and discrete distributions for forecasted production of stochastic units, forecasted consumption as well as the corresponding forecast errors, and represent the gradual realization of uncertainty by a so-called scenario tree \cite{birge2011}. Furthermore, by linearization of the non-linear terms in the formulations and using the optimality conditions in the bi-level problem, we cast the proposed expansion planning models as  mixed-integer linear programs.

In summary, the main contributions of this paper are:
\begin{enumerate}
\item The proposal of a stochastic programming framework to determine the optimal power system expansion plan that takes into account both the variability of the stochastic production and the corresponding forecast errors.
\item The analysis of how such forecast errors impact the optimal expansion decisions.
\item The investigation of how market design influences the expansion planning of a power system with high penetration of stochastic production.
\end{enumerate}

The paper is structured as follows. Section \ref{Intro} introduces and motivates the expansion planning problem in power systems with high penetration of renewable energy. In Section \ref{ModForm} we first present the expansion planning problem assuming perfect forecasts of stochastic production. We then extend this problem to include forecast errors under the two  market designs outlined above. Section \ref{IllEx} provides an illustrative four-node example of the expansion problem. A more realistic 24-node case study is presented in Section \ref{CaseStudy}. Finally, Section \ref{Con} concludes the paper.

\section{Modeling and formulation}\label{ModForm}

In this paper, we consider an existing power system that consists of a set $\mathcal{G}$ of generating units, a set $\mathcal{L}$ of loads and a set $\mathcal{F}$ of transmission lines. The proposed stochastic programming models determine the optimal decisions among a set $\mathcal{P}$ of available capacity expansion projects to minimize the sum of operating and investment costs while ensuring that at least $\eta\%$ of the electricity consumption is covered by renewable electricity production. In order to model different states of the power system throughout the decision horizon, a set of scenarios $s$ is considered. The probability of each scenario is denoted by $\pi_s$, such that $\sum_s \pi_s = 1$.

Expansion projects include both dispatchable and stochastic generating units as well as new transmission lines. Each expansion project $i\in\mathcal{P}$ is limited by a maximum capacity $x^{\rm max}_i$ and involves a fixed investment cost $q^0_i$ and a variable investment cost $q_i$. Investment decisions are modeled by a binary variable $u_i$, that is equal to 1 if the project is carried out and 0 otherwise, and by the optimal capacity of the project $x_i$. The total investment cost of the system is computed as $\sum_{i\in\mathcal{P}}(q^0_{i}u_{i}\Plus q_{i}x_{i})$.

Generating units $i\in\mathcal{G}$ are modeled by a maximum capacity $x_i$, a marginal cost $c_i$ and a capacity factor $\widehat{\rho}_{is}$ that depends on the realized scenario $s$. Note that $x_i$ is a known and fixed parameter for existing units, while it represents a variable for new generating units. Parameter $\widehat{\rho}_{is}$ can be used to characterize both the occurrence of unexpected unit failures and the variable generation from stochastic generating units, such as wind or solar power plants.

Consumption units $i\in\mathcal{L}$ are characterized by a known peak load $x_i$, a capacity factor $\widehat{\rho}_{is}$ that can be used to compute the load level for each scenario $s$, and a utility $c_i$. Observe that if loads are assumed to be inflexible, then $c_i=\Minus v^{LS}_i \: \forall i\in\mathcal{L}$, where $v^{LS}_i$ stands for the cost of involuntary load shedding. Otherwise, $c_i$ may represent the flexibility of the loads according to their utility function.

The transmission network is modeled using DC power flow equations, in which transmission lines have a maximum capacity $x_i$ and a susceptance $b_i$. Like for generating units, $x_i$ is a known and fixed parameter for the existing lines. The capacity factor $\widehat{\rho}_{is}$ may model line failures. The parameter $a_{in}$ denotes the line-bus indicator, and $\widehat{\delta}_{ns}$ stands for the voltage angle at bus $n$ and scenario $s$.

According to these data and following a network-constrained economic dispatch, the market operator determines the production of generating units ($\widehat{p}_{is} \: \forall i\in\mathcal{G}$), the consumption of loads ($\widehat{p}_{is} \: \forall i\in\mathcal{L}$) and the power flow through the transmission lines ($\widehat{p}_{is} \: \forall i\in\mathcal{F}$) that minimize the operating cost for each scenario realization $s$.

Due to the high predictability of electricity demand and the low penetration of stochastic production, classical expansion planning models systematically disregard the effect of forecast errors. In this vein, we present a generic expansion planning formulation that accounts for the variability of demand and stochastic production throughout the planning horizon, but ignores the associated forecast errors. 
\begin{subequations} \label{ExpOnlyDA} \begin{align}
& {\rm Minimize}_{u_{i},x_{i},\widehat{p}_{is},\widehat{\delta}_{ns}} \nonumber \\
& \qquad \qquad  \sum_{i\in\mathcal{P}} \left(q^0_{i}u_{i}
\Plus q_{i}x_{i}\right) \Plus \hspace{-3mm} \sum_{i\in\mathcal{G} \cup \mathcal{L},s} \hspace{-2mm} T\pi_{s}c_i\widehat{p}_{is} \label{ExpOnlyDA_TotalCost}\displaybreak[1] \\
& {\rm subject \: to} \hspace{8cm} \nonumber \\
& \sum_{i\in \mathcal{W},s } \hspace{-1mm} \pi_s\widehat{p}_{is} \geq \eta\sum_{i\in \mathcal{L},s} \pi_s\widehat{p}_{is} \label{ExpOnlyDA_MinWindShare} \displaybreak[1]\\
& 0 \leq x_{i} \leq u_{i}x^{\rm max}_{i}, \quad  \forall  i\in \mathcal{P} \label{ExpOnlyDA_MaxCap}\displaybreak[1]\\
& 0 \leq \widehat{p}_{is} \leq x_{i}\widehat{\rho}_{is}, \quad  \forall i\in\mathcal{G} \cup \mathcal{L}, \forall s \label{ExpOnlyDA_MaxP_DA}\displaybreak[1]\\
& \widehat{p}_{is} = u_i b_i \sum_n a_{in}\widehat{\delta}_{ns}, \quad \forall i\in\mathcal{F}, \forall s \label{ExpOnlyDA_Flow_DA}\displaybreak[1]\\
& \Minus x_{i}\widehat{\rho}_{is}  \leq \widehat{p}_{is}  \leq x_{i}\widehat{\rho}_{is}, \quad \forall i\in\mathcal{F}, \forall s   \label{ExpOnlyDA_MaxFlow_DA}\displaybreak[1]\\
& \hspace{-1mm} \sum_{i\in\mathcal{G}_n} \widehat{p}_{is} = \hspace{-1mm} \sum_{i\in\mathcal{L}_n} \widehat{p}_{is} \Plus \hspace{-1mm} \sum_{i\in\mathcal{F}} \hspace{0mm} a_{in}\widehat{p}_{is}, \quad \forall n, \forall s \label{ExpOnlyDA_Balance_DA}\displaybreak[1]\\
& \widehat{\delta}_{n_1s}=0, \quad \forall s \label{ExpOnlyDA_SlackBus_DA} \displaybreak[1] \\
& u_i\in\{0,1\}, \quad \forall i\in\mathcal{P} \label{ExpOnlyDA_Binary}\displaybreak[1]\\
& u_i = 1, \quad \forall i \notin   \mathcal{P} \label{ExpOnlyDA_Existing}\displaybreak[1]
\end{align}
\end{subequations}

Objective function \eqref{ExpOnlyDA_TotalCost} minimizes current investment cost plus future operating costs throughout the planning horizon. Constraint \eqref{ExpOnlyDA_MinWindShare} ensures that at least $\eta\%$ of the demand is covered with electricity produced by renewable generating units $\mathcal{W}$. The capacity of new generating units and transmission lines is bounded by the maximum capacity of each project in $\mathcal{P}$ through \eqref{ExpOnlyDA_MaxCap}. The dispatch of generating units and loads is limited by their corresponding capacities in \eqref{ExpOnlyDA_MaxP_DA}. Likewise, the power flow defined in \eqref{ExpOnlyDA_Flow_DA} is bounded by the capacities of the transmission lines $\mathcal{F}$ through \eqref{ExpOnlyDA_MaxFlow_DA}. Constraint \eqref{ExpOnlyDA_Balance_DA} ensures the power balance at each node. Equation \eqref{ExpOnlyDA_SlackBus_DA} arbitrarily sets the voltage angle at bus $n_1$ to 0. Finally, \eqref{ExpOnlyDA_Binary} and \eqref{ExpOnlyDA_Existing} are binary variable declarations.

Optimization model \eqref{ExpOnlyDA} is a two-stage stochastic programming problem where first-stage variables are the investment decisions ($u_i,x_i$), the uncertain parameters are the capacity factors throughout the planning horizon ($\widehat{\rho}_{is}$), and the second-stage variables are dispatch decisions ($\widehat{p}_{is},\widehat{\delta}_{is}$), which depend on each particular realization of the capacity factors. After linearizing the product of continuous and binary variables according to the procedure presented in \cite{floudas1995}, this model can be cast as a mixed-integer linear programming problem and thus solved using commercial software.

This model of the electricity market could represent two situations according to how power systems are operated today. One could think of it as a real-time market in which all generating units are assumed to be completely flexible and able to instantaneously adapt their output to the status of the system. Another interpretation could be a day-ahead market cleared with perfect forecasts of the demand and stochastic production. However, most generating units have technical constraints regarding their response time to unexpected events and the eventual values of demand and stochastic production differ from the forecasted ones.

In this paper, we propose to overcome these shortcomings by including an additional decision stage into model \eqref{ExpOnlyDA}. This allows us to more accurately characterize the stochastic parameters involved by taking into account probability distributions of the corresponding forecast errors, and to model technical constraints on the flexibility provided by generating units. The proposed modeling of the functioning of the market in two stages is to be interpreted as the clearing of a day-ahead market taking place 24-36 hours in advance according to forecasts followed by a balancing market that copes with forecast errors. In balancing markets, generating units can submit quantity-price offers on how much they are willing to deviate with respect to their day-ahead dispatch. For example, a generating unit $i\in\mathcal{G}$ can offer up- and down-balancing energy up to $x^{\Plus}_ix_i$ and $x^{\Minus}_ix_i$ MW at a cost of $c^{\Plus}_i$ and $c^{\Minus}_i$ \$$/$MWh, respectively. We can also assume that inflexible loads can work as up-balancing resources through curtailment. Likewise, stochastic generating units can provide down-balancing service through spillage.

Considering both a day-ahead and a balancing market raises further questions on whether the level of coordination between these two markets may influence the expansion planning of a power system. To address this issue, we propose two different  optimization models that differ in whether forecast errors are accounted for or not when deciding on the day-ahead dispatch. First, and following the proposal described in \cite{Pritchard2010,Morales2013}, we present an expansion planning model in which the functioning of the day-ahead and the balancing markets is co-optimized:
\begin{subequations} \label{ExpStoc} \begin{align}
& {\rm Minimize}_{u_{i},x_{i},\widehat{p}_{is},\widehat{\delta}_{ns}, p^{\Plus}_{isr},p^{\Minus}_{isr},\widetilde{\delta}_{nsr}} \nonumber \\
& \hspace{-2mm} \sum_{i\in\mathcal{P}} \left(q^0_{i}u_{i}
\Plus q_{i}x_{i}\right) \Plus \hspace{-2mm} \sum_{i\in\mathcal{G} \cup \mathcal{L},s} \hspace{-2mm} T\pi_{s} \left( c_i\widehat{p}_{is} \Plus \sum_{r}  \pi_{sr} \left(c^{\Plus}_ip^{\Plus}_{isr} \Minus c^{\Minus}_ip^{\Minus}_{isr} \right) \right) \label{ExpStoc_TotalCost}\displaybreak[1] \\
& {\rm subject \: to} \nonumber \\
& \sum_{i\in \mathcal{W},sr } \hspace{-1mm} \pi_s \pi_{sr} \widetilde{p}_{isr}  \geq \eta\sum_{i\in \mathcal{L},sr} \pi_s \pi_{sr} \widetilde{p}_{isr}  \label{ExpStoc_MinWindShare} \displaybreak[1]\\
& \eqref{ExpOnlyDA_MaxCap}-\eqref{ExpOnlyDA_Existing} \displaybreak[1]\\
& 0 \leq \widetilde{p}_{isr} \leq x_{i} \widetilde{\rho}_{isr}, \quad \forall i\in\mathcal{G} \cup \mathcal{L},\forall s, \forall r \label{ExpStoc_MaxP_B}\displaybreak[1]\\
& 0 \leq p^{\Plus}_{isr} \leq x^{\Plus}_ix_i, \quad \forall i\in\mathcal{G} \cup \mathcal{L},\forall s, \forall r \label{ExpStoc_MaxPUp}\displaybreak[1]\\
& 0 \leq p^{\Minus}_{isr} \leq x^{\Minus}_ix_i, \quad \forall i\in\mathcal{G} \cup \mathcal{L},\forall s, \forall r \label{ExpStoc_MaxPDo}\displaybreak[1]\\
& \widetilde{p}_{isr} = u_i b_i \sum_n a_{in}\widetilde{\delta}_{nsr}, \quad \forall i\in\mathcal{F}, \forall s, \forall r \label{ExpStoc_Flow_B}\displaybreak[1]\\
& \Minus x_{i}\widetilde{\rho}_{isr}  \leq \widetilde{p}_{isr}  \leq x_{i}\widetilde{\rho}_{isr}, \quad \forall i\in\mathcal{F}, \forall s, \forall r   \label{ExpStoc_MaxFlow_B}\displaybreak[1]\\
& \hspace{-1mm} \sum_{i\in\mathcal{G}_n} \widetilde{p}_{isr} = \hspace{-1mm} \sum_{i\in\mathcal{L}_n} \widetilde{p}_{isr} \Plus \hspace{-1mm} \sum_{i\in\mathcal{F}} \hspace{0mm} a_{in}\widetilde{p}_{isr}, \quad \forall n, \forall s, \forall r \label{ExpStoc_Balance_B}\displaybreak[1]\\
& \widetilde{\delta}_{n_1sr}=0, \quad \forall s, \forall r, \label{ExpStoc_SlackBus_B} \displaybreak[1]
\end{align}
\end{subequations}

\noindent where the final dispatch $\widetilde{p}_{isr} \forall i\in\mathcal{G} \cup \mathcal{L} $ is equal to the day-ahead dispatch $\widehat{p}_{is}$ plus the deployed up-balacing power $p^{\Plus}_{isr}$ minus the down-balancing power $p^{\Minus}_{isr}$, i.e., $\widetilde{p}_{isr}=\widehat{p}_{is} \Plus p^{\Plus}_{isr} \Minus p^{\Minus}_{isr}$. Likewise, $\widetilde{p}_{isr} \forall i\in\mathcal{F}$ represents the power flow through transmission line $i$ in balancing scenario $sr$.

Objective function \eqref{ExpStoc_TotalCost} minimizes the investment cost plus the expected operation cost, which includes both the day-ahead dispatch and the balancing redispatch costs. Equation \eqref{ExpStoc_MinWindShare} ensures the minimum penetration of renewable electricity production. Constraints \eqref{ExpOnlyDA_MaxCap}-\eqref{ExpOnlyDA_Existing} in the simplified expansion model \eqref{ExpOnlyDA} are also needed here to model the functioning of the day-ahead market. In the same fashion as in formulation \eqref{ExpOnlyDA}, equations \eqref{ExpStoc_MaxP_B}, \eqref{ExpStoc_MaxPUp}, \eqref{ExpStoc_MaxPDo}, \eqref{ExpStoc_Flow_B}, \eqref{ExpStoc_MaxFlow_B}, \eqref{ExpStoc_Balance_B}, and \eqref{ExpStoc_SlackBus_B} limit the dispatch and re-dispatch of generating units, compute the power flow through the transmission lines, impose bounds on the power flows, ensure the power balance at each bus, and arbitrarily set bus $n_1$ as the reference node at the balancing stage, respectively.

Optimization problem \eqref{ExpStoc} is a three-stage stochastic optimization problem in which the uncertainty reveals over time as follows: at the day-ahead stage, only the forecast of the capacity factors is issued; and at the balancing stage we assume that the actual values of such parameters are fully known. Variables are therefore divided into three types: the first-stage variables are the investment decisions ($u_i,x_i$), which are determined facing all the sources of uncertainty involved; the second-stage variables are the dispatch decisions at the day-ahead stage ($\widehat{p}_{is},\widehat{\delta}_{ns}$) which are made conditional on the forecast values of the capacity factors ($\widehat{\rho}_{is}$); and the third-stage variables are the adjustments to the day-ahead dispatch  ($p^{\Plus}_{isr},p^{\Minus}_{isr},\widetilde{\delta}_{nsr}$), which depend on the realization of the actual capacity factors ($\widetilde{\rho}_{isr}$). By linearizing the product of continuous and binary variables, this model can likewise be solved as a mixed-integer linear problem.

Note that model \eqref{ExpStoc} decides on the optimal expansion plan assuming an ideal market in which forecast errors of stochastic units are handled as efficiently as possible. In order to investigate the extent to which the expansion plan is affected by how forecast errors are processed by the market, we build the alternative expansion planning model \eqref{ExpConv}. Contrary to model \eqref{ExpStoc}, model \eqref{ExpConv} disregards the potential impact of forecast errors on the balancing costs, thus failing to provide a day-ahead dispatch that makes an efficient use of the available flexible generation \cite{Morales2013}.
\begin{subequations} \label{ExpConv}
\begin{align}
& {\rm Minimize}_{u_{i},x_{i},\widehat{p}_{is},\widehat{\delta}_{ns}, p^{\Plus}_{isr},p^{\Minus}_{isr},\widetilde{\delta}_{nsr}} \nonumber \\
& \hspace{-2mm} \sum_{i\in\mathcal{P}} \left(q^0_{i}u_{i}
\Plus q_{i}x_{i}\right) \Plus \hspace{-2mm} \sum_{i\in\mathcal{G} \cup \mathcal{L},s} \hspace{-2mm} T\pi_{s} \left( c_i\widehat{p}_{is} \Plus \sum_{r}  \pi_{sr} \left(c^{\Plus}_ip^{\Plus}_{isr} \Minus c^{\Minus}_ip^{\Minus}_{isr} \right) \right)  \label{ExpConv_TotalCost}\displaybreak[1] \\
& {\rm subject \: to} \hspace{8cm} \nonumber \\
& \eqref{ExpStoc_MinWindShare}-\eqref{ExpStoc_SlackBus_B}
\end{align}
\begin{empheq}[left={\begin{matrix} \widehat{p}_{is} \\ \widehat{\delta}_{ns} \end{matrix} \hspace{-1mm}\in\hspace{-1mm} {\rm arg}}\empheqlbrace,right=\empheqrbrace \forall s]{align}
& {{\rm  Minimize}_{\widehat{p}_{is},\widehat{\delta}_{ns}}} \sum_{i\in\mathcal{G} \cup \mathcal{L},s} \hspace{-2mm} c_i\widehat{p}_{is}  \label{ExpConv_CostDA} \displaybreak[1]\\
& \textrm{subject to} \nonumber \displaybreak[1]\\
& 0 \leq \widehat{p}_{is} \leq x_{i}\widehat{\rho}_{is}:\underline{\alpha}_{is},\overline{\alpha}_{is}, \forall i\in\mathcal{G} \cup \mathcal{L} \label{ExpConv_MaxP_DA}\displaybreak[1]\\
& \widehat{p}_{is} = u_i b_i \sum_n a_{in}\widehat{\delta}_{ns}:\phi_{is}, \forall i\in\mathcal{F} \label{ExpConv_Flow_DA}\displaybreak[1]\\
& \Minus x_{i}\widehat{\rho}_{is}  \leq \widehat{p}_{is}  \leq x_{i}\widehat{\rho}_{is}:\underline{\theta}_{is},\overline{\theta}_{is}, \forall i\in\mathcal{F}   \label{ExpConv_MaxFlow_DA}\displaybreak[1]\\
& \hspace{-1mm} \sum_{i\in\mathcal{G}_n} \widehat{p}_{is} = \hspace{-1mm} \sum_{i\in\mathcal{L}_n} \widehat{p}_{is} \Plus \hspace{-1mm} \sum_{i\in\mathcal{F}} \hspace{0mm} a_{in}\widehat{p}_{is}:\hspace{-1mm}\lambda_{ns}, \forall n \label{ExpConv_Balance_DA}\displaybreak[1]\\
& \widehat{\delta}_{n_1s}=0: \xi_s \label{ExpConv_SlackBus_DA} \displaybreak[1]
\end{empheq}
\end{subequations}

Observe that the objective function and all constraints of problem \eqref{ExpStoc} are included in \eqref{ExpConv}. However, a new set of constraints \eqref{ExpConv_CostDA}-\eqref{ExpConv_SlackBus_DA} is added to this formulation in order to impose that the day-ahead decisions $\widehat{p}_{is},\widehat{\delta}_{ns}$ are those that  minimize the day-ahead cost alone, with no account taken of the potentital impact of these decisions on the ensuing balancing operation of the power system. The lower-level optimization problem represents the clearing of the day-ahead market as formulated in \eqref{ExpOnlyDA}, but also includes the dual variables corresponding to each constraint after a colon. Therefore, problem \eqref{ExpConv} is a more constrained version of problem \eqref{ExpStoc}. Problem \eqref{ExpConv} has a bilevel structure, and to solve it, the lower level problem \eqref{ExpConv_CostDA}-\eqref{ExpConv_SlackBus_DA} is replaced with its KKT conditions or, alternatively, with its primal constraints, dual constraints, plus the strong duality condition  as follows \cite{motto2005}:
\begin{subequations} \label{ExpConv2} \begin{align}
& {\rm Minimize}_{u_{i},x_{i},\widehat{p}_{is},\widehat{\delta}_{ns}, p^{\Plus}_{isr},p^{\Minus}_{isr},\widetilde{\delta}_{nsr}} \nonumber \\
& \hspace{-2mm} \sum_{i\in\mathcal{P}} \left(q^0_{i}u_{i}
\Plus q_{i}x_{i}\right) \Plus \hspace{-2mm} \sum_{i\in\mathcal{G} \cup \mathcal{L},s} \hspace{-2mm} T\pi_{s} \left( c_i\widehat{p}_{is} \Plus \sum_{r}  \pi_{sr} \left(c^{\Plus}_ip^{\Plus}_{isr} \Minus c^{\Minus}_ip^{\Minus}_{isr} \right) \right) \label{ExpConv2_TotalCost}\displaybreak[1] \\
& {\rm subject \: to} \hspace{8cm} \nonumber \\
& \eqref{ExpStoc_MinWindShare}-\eqref{ExpStoc_SlackBus_B} \displaybreak[1] \\
& \underline{\alpha}_{is} \Plus\: \overline{\alpha}_{is} \Plus\: \lambda_{n_{i}s} = c_i, \quad \forall i\in\mathcal{G}, \forall s \label{ExpConv2_dual1} \displaybreak[1]\\
& \underline{\alpha}_{is} \Plus\: \overline{\alpha}_{is} \Minus\: \lambda_{n_{i}s} = c_i, \quad \forall i\in\mathcal{L}, \forall s \label{ExpConv2_dual2} \displaybreak[1]\\
& \phi_{is} \Plus\: \underline{\theta}_{is} \Plus\: \overline{\theta}_{is} \Minus \sum_n a_{in}\lambda_{ns} = 0, \quad \forall i\in\mathcal{F}, \forall s \label{ExpConv2_dual3} \displaybreak[1] \\
& \Minus \sum_{i\in\mathcal{F}} u_ib_ia_{in}\phi_{is} \Plus \left( \xi_s\right)_{n=n_1} = 0, \quad \forall n, \forall s \label{ExpConv2_dual4} \displaybreak[1] \\
& \sum_{i\in\mathcal{G} \cup \mathcal{L},s} \hspace{-2mm} c_i\widehat{p}_{is} = \hspace{-3mm} \sum_{i\in\mathcal{G} \cup \mathcal{L},s} \hspace{-2mm} x_i\widehat{\rho}_{is} \overline{\alpha}_{is} \Plus \sum_{i\in\mathcal{F}} x_i\widehat{\rho}_{is} \left(\overline{\theta}_{is} \Minus \underline{\theta}_{is}  \right), \forall s, \label{ExpConv2_dual5} \displaybreak[1]
\end{align} \end{subequations}

\noindent where $n_i$ denotes the bus to which device $i\in\mathcal{G} \cup \mathcal{L}$ is connected. Equations \eqref{ExpConv2_dual1}-\eqref{ExpConv2_dual4} are the dual constraints corresponding to the lower-level problem \eqref{ExpConv_CostDA}-\eqref{ExpConv_SlackBus_DA}, and constraint \eqref{ExpConv2_dual5}  formulates the strong primal-dual condition. Note that equation \eqref{ExpConv2_dual5} includes products of two continuous variables, namely, the capacity of the projects $x_i$ and some dual variables $\overline{\alpha}_{is},\overline{\theta}_{is}, \underline{\theta}_{is}$. In order to linearize these terms, we use the binary expansion of the capacity investment decision $x_i$ as follows:
\begin{equation}
x_i = \sum_{b=1}^{N^B_i} v_{ib}S_i2^{(b-1)}, \quad \forall i\in\mathcal{P},
\end{equation}

\noindent where $v_{ib}$ is a binary variable for each block $b$, $S_i$ is the block size (which is a parameter) and the number of blocks for each project is computed as

\begin{equation}
N^B_{i} = floor\left(log_2\left(\dfrac{x^{\rm max}_{i}}{S_{i}} \right)\right) + 1.
\end{equation}

\noindent In doing so, problem \eqref{ExpConv2} is also formulated as a mixed-integer linear optimization problem.

\section{Uncertainty characterization}\label{UncerCharact}

Uncertainty comes into expansion models~\eqref{ExpOnlyDA}, \eqref{ExpStoc}, and \eqref{ExpConv} through the capacity factors $\widehat{\rho}_{is}$ and $\widetilde{\rho}_{isr}$. In this section, we briefly describe how scenarios for these capacity factors are generated. For conciseness, though, we limit ourselves to the case of wind power uncertainty, as the procedure to generate scenarios for load and equipment failures runs in a similar manner. Furthermore, the scenario modeling approach used in this paper is analogous to the one described in \cite{Pineda2014}, to which the reader is hereby referred.

The scenario generation procedure proceeds in two steps:
\begin{enumerate}
  \item First, a number $N_{\text s}$ of per-unit power values ($\widehat{\rho}_{is}$), each representing a 24-hour ahead predicted power output of wind location $i\in\mathcal{W}$ (in per unit), are sampled from the stationary probability distribution that characterizes the per-unit power production from this wind site, and that can be obtained using historical data \cite{Holttinen2005}.
  \item Second, for each of these per-unit power values $\widehat{\rho}_{is}$, a set of $N_{\text r}$ prediction errors ($\widetilde{\rho}_{isr}$) are sampled from the per-unit forecast error distribution that also characterizes wind location $i\in\mathcal{W}$. Here, we assume that wind power forecast errors follow a Beta distribution, as in \cite{Fabbri2005}. Furthermore, the characteristic parameters of the Beta distribution to be used are functions of the predicted wind power output (i.e., $\widehat{\rho}_{is}$) and the length of the forecast horizon. These functions are also provided in \cite{Fabbri2005}.
\end{enumerate}

As stated above, this two-step procedure also applies to load uncertainty and equipment failures, by just considering the appropriate probability distributions. Note also that scenario reduction techniques can be applied to reduce the size of the scenario set while maintaining most of the statistical features of the stochastic parameters \cite{Morales2009}.

\section{Illustrative example}\label{IllEx}

Next a small example is used to provide intuition about the three expansion planning models previously described. Despite its reduced size, this example provides evidence to the fact that the limited predictability of renewable energy sources impacts the optimal expansion planning of a power system. Furthermore, the magnitude of this impact is contingent on the design of the electricity market that governs the short-term operation of the system.

\subsection{Data}

\begin{figure}[h!]
    \centering  \includegraphics[scale=0.6]{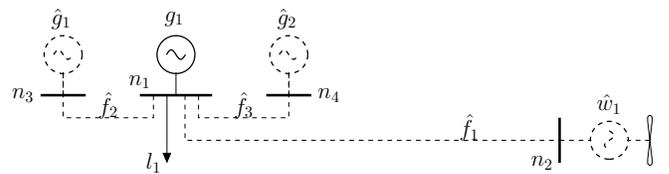}
    \caption{4-bus power system}
    \label{fig:4bus}
\end{figure}

\begin{figure}[h!]
\centering 
\subfigure[Demand]{\includegraphics[scale=0.27]
{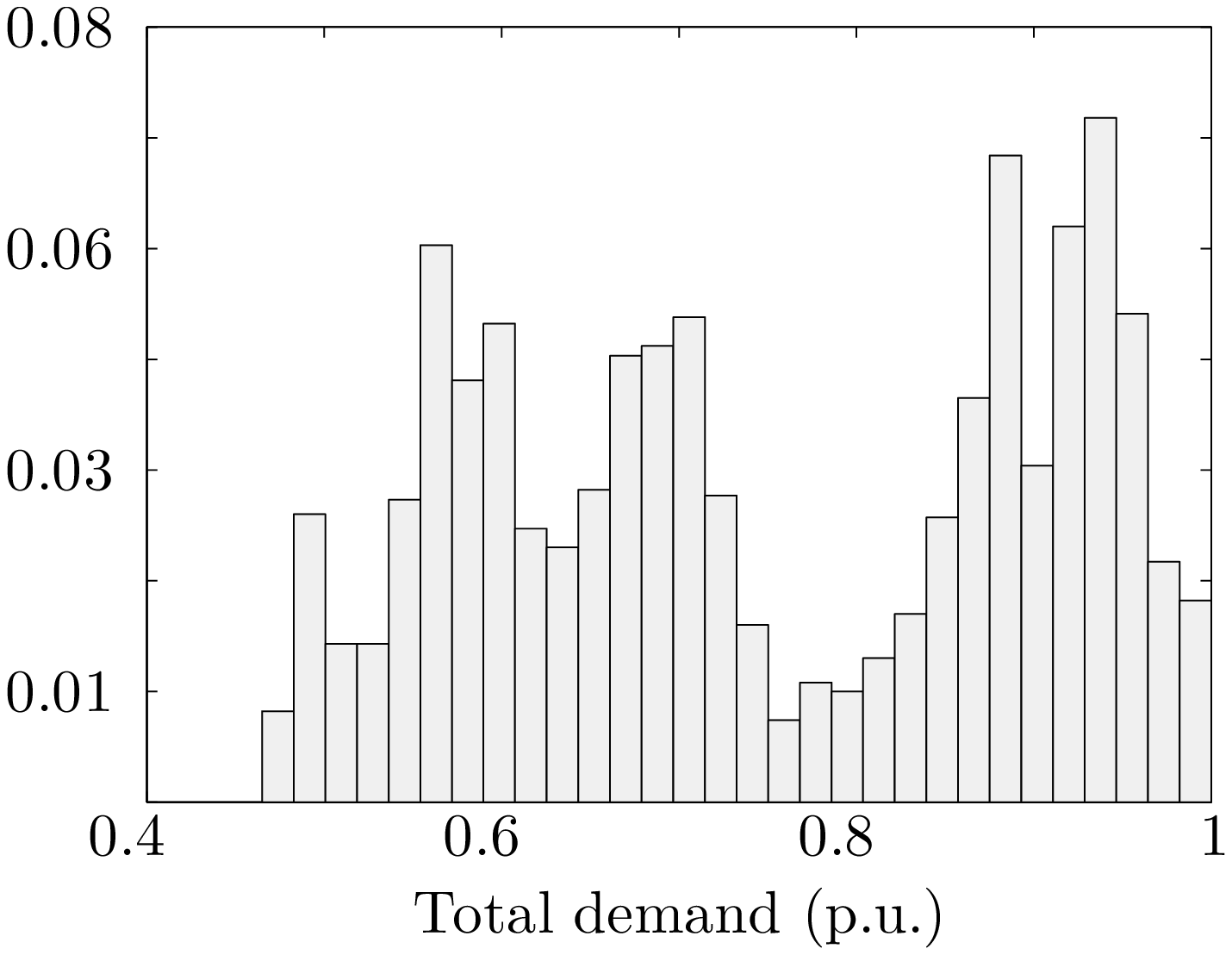}\label{fig:load}} 
\subfigure[Wind]{\includegraphics[scale=0.27]
{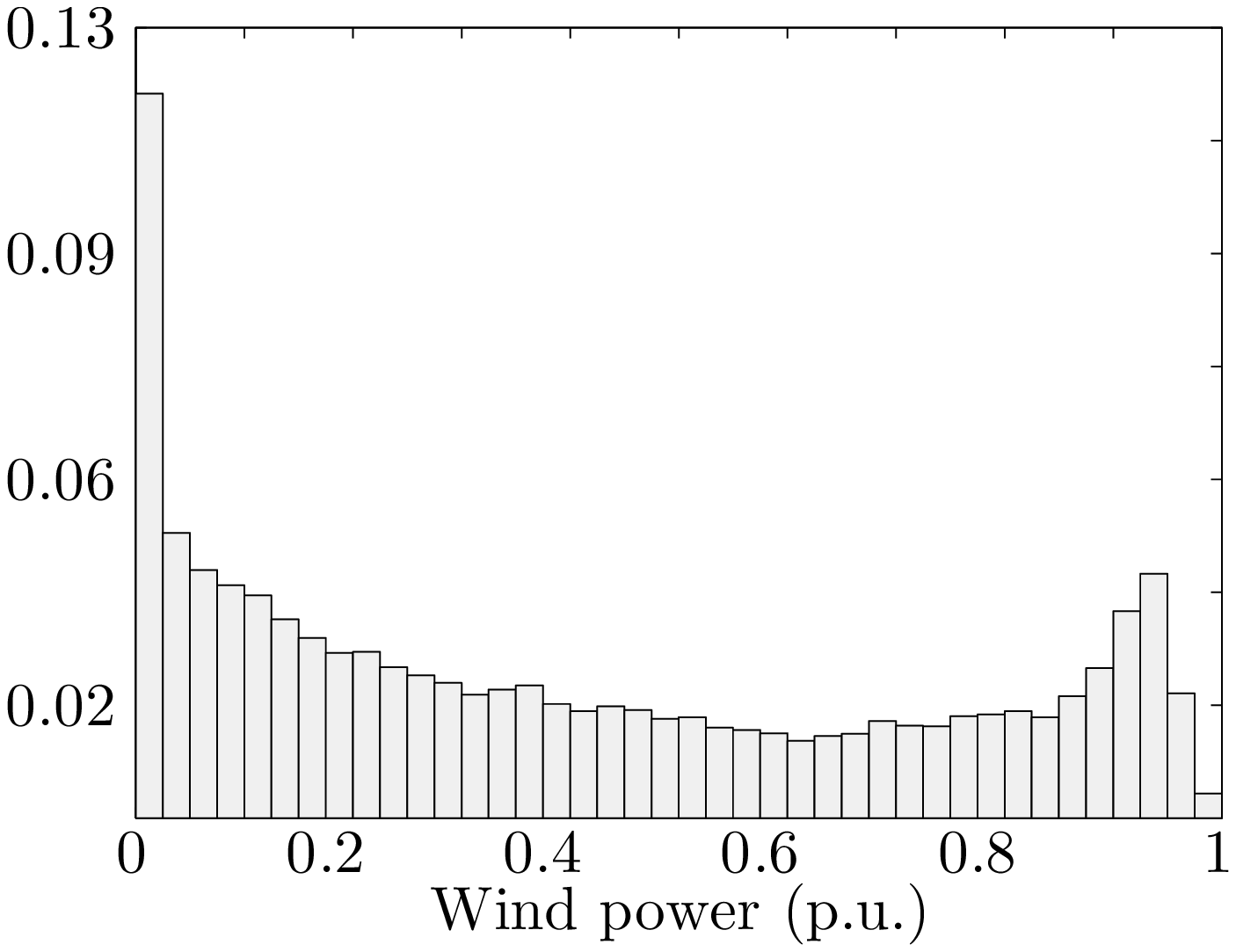}\label{fig:wind}} 
\caption{Probability distributions}\end{figure}

Fig. \ref{fig:4bus} shows a small power system that initially consists of one load $l_{1}$ and an \emph{inflexible} generating unit $g_1$ (solid lines). Fig. \ref{fig:load} plots the probability distribution of the per-unit consumption at the pre-existing node $n_1$, with a peak load of 500 MW. The marginal cost and capacity of $g_1$ are equal to \$10$/$MWh and 500 MW, respectively.

In order to increase the contribution of renewables to the electricity supply (e.g., to comply with a certain target), the system in Fig.~\ref{fig:4bus} can be expanded through the following available projects (dashed lines):

\begin{itemize}
\item[-] The installation of a wind farm with a capacity of up to 1000 MW at $n_2$. This site is characterized by the per-unit wind power production depicted in Fig. \ref{fig:wind}. The fixed and variable annualized investment costs associated with this project are $q^0_i=\$25\,000$ and $q_i=\$50\,000/{\rm MW}$.
\item[-] Two new thermal generating units, $\hat{g}_1$ and $\hat{g}_2$, at nodes $n_3$ and $n_4$, respectively. The characteristics of these units are collated in Table \ref{table:2busgendata}, with the units of the different parameters provided in the Notation section. Note that both units have the same marginal cost $c_i$ and the same investment costs $q^0_i,q_i$, and that, unlike the pre-existing generator $g_1$, they are \emph{flexible}, i.e., they both provide energy for balancing. More specifically, unit $\hat{g}_2$ offers more expensive upward regulation than unit $\hat{g}_1$, but is also willing to repurchase energy in the balancing market at a much higher price.
\item[-] Three new transmission lines, $\hat{f}_1$, $\hat{f}_2$ and $\hat{f}_3$, to connect the pre-existing bus $n_1$ with the wind farm, unit $\hat{g}_1$, and unit $\hat{g}_2$, respectively. The data pertaining to the new transmission lines is provided in Table \ref{table:2buslinedata}. Lines $\hat{f}_2$ and $\hat{f}_3$ have the same length, which is in turn five times shorter than the length of $\hat{f}_1$. Correspondingly, the investment costs of $\hat{f}_2$ and $\hat{f}_3$ are five times lower than those of $\hat{f}_1$.
\end{itemize}

\begin{table}[htb] \centering \renewcommand*{\arraystretch}{1.5}
\caption{New generating units (example)}
\begin{tabular}{c c c c c c c c c}
  \hline
  {\vphantom{\Large{$A_p$}}} & $x^{max}_i$ & $c_{i}$  & $x^{\Plus}_i$ & $c^{\Plus}_i$ & $x^{\Minus}_i$ & $c^{\Minus}_i$ & $q^0_{i}$ & $q_{i}$ \\
  \hline
  $\hat{g}_1$ & 250 & 20   & 100  & 21   & 100  & 0   & 20000 & 25000\\
  $\hat{g}_2$ & 250 & 20   & 100  & 22   & 100  & 20  & 20000 & 25000\\
  \hline
\end{tabular}\label{table:2busgendata}\end{table}

\begin{table}[htb] \centering \renewcommand*{\arraystretch}{1.5}
\caption{New transmission lines (example)}
\begin{tabular}{c c c c c}
  \hline
  {\vphantom{\Large{$A_p$}}} & $x^{max}_i$ & $b_{i}$  &  $q^0_{i}$ & $q_{i}$ \\
  \hline
  $\hat{f}_1$ & 500 & 20   & 25625 & 90\\
  $\hat{f}_2$ & 100 & 20   & 5125  & 18\\
  $\hat{f}_3$ & 100 & 20   & 5125  & 18\\
  \hline
\end{tabular}\label{table:2buslinedata}\end{table}

According to the notation used in Section~\ref{ModForm}, we have the sets $\mathcal{F}=\{\hat{f}_1,\hat{f}_2,\hat{f}_3\}$, $\mathcal{G}=\{g_1,\hat{g}_1,\hat{g}_2,\hat{w}_1 \}$, $\mathcal{L}=\{l_1 \}$, $\mathcal{P}=\{\hat{f}_1,\hat{f}_2,\hat{f}_3,\hat{g}_1,\hat{g}_2,\hat{w}_1 \}$, $\mathcal{W}=\{\hat{w}_1 \}$.

For simplicity, unexpected failures of thermal generating units and transmission lines are disregarded here, i.e., $\widehat{\rho}_{is}=\widetilde{\rho}_{isr}=1, \: \forall i\in\mathcal{F} \cup \{ \mathcal{G} \backslash \mathcal{W} \}, \forall s, \forall r$. Moreover, load $l_1$ is assumed to be inelastic with $V^{LS} = \$500/$MWh, and demand forecast errors are disregarded. Finally, the capacity of the available expansion projects is discretized in 1-MW blocks.

The variability of the forecasts for wind power production and demand throughout the planning horizon is approximated by 20 scenarios. Likewise, conditional on each scenario at the day-ahead stage, a set of 30 scenarios is generated to characterize the forecast errors of the wind power production as explained in Section \ref{UncerCharact}. Therefore, a final set of 600 scenarios is considered in this example.

\subsection{Optimal expansion planning results}

The plans to expand the power system, together with their associated total investment costs from expansion planning models~\eqref{ExpOnlyDA}, \eqref{ExpStoc}, and \eqref{ExpConv} are compared in Table \ref{table:resultsexpansion} as a function of a minimum target of wind power penetration. Recall that model~\eqref{ExpOnlyDA} disregards the forecast errors of wind power production, while these are considered in both \eqref{ExpStoc} and \eqref{ExpConv}. However, whereas \eqref{ExpStoc} assumes that forecast errors are efficiently handled by a market organization that optimally coordinates the day-ahead generation scheduling with the subsequent balancing operation of the power system,
\eqref{ExpConv} presupposes no coordination at all. Note also that the optimal investments in lines $\hat{f}_2$ and $\hat{f}_3$ are not included here because they turn out to be in all cases equal to the capacities of generating units $\hat{g}_1$ and $\hat{g}_2$, respectively. 

\begin{table}[htb] \centering \renewcommand*{\arraystretch}{1.5}
\caption{Optimal expansion plans (example). Capacities in MW and investment cost in m\$.}
\begin{tabular}{c c c c c c c}
  \hline
  && \multicolumn{5}{c}{Renewable target}\\
  && 10\% & 20\% & 30\% & 40\% & 50\% \\
  \hline
  & model \eqref{ExpOnlyDA}& 98  & 196 & 294 & 392 & 557 \\  
  Cap. $\hat{w}_1$ & model \eqref{ExpStoc}  & 98  & 196 & 294 & 413 & 580 \\  
  & model \eqref{ExpConv} & 104 & 208 & 312 & 441 & 629 \\  
  \hline
  & model \eqref{ExpOnlyDA} & 98  & 196 & 294 & 392 & 456 \\  
  Cap. $\hat{f}_1$ & model \eqref{ExpStoc} & 97  & 193 & 290 & 369 & 456 \\  
  & model \eqref{ExpConv} & 94  & 188 & 282 & 377 & 454 \\ 
  \hline
  & model \eqref{ExpOnlyDA} & 0   & 0   & 0   & 0   & 0   \\  
  Cap. $\hat{g}_1$ & model \eqref{ExpStoc} & 0   & 0   & 0   & 0   & 0   \\   
  & model \eqref{ExpConv} & 54  & 107 & 161 & 195 & 199 \\ 
  \hline
  & model \eqref{ExpOnlyDA} & 0   & 0   & 0   & 0   & 0   \\   
  Cap. $\hat{g}_2$ & model \eqref{ExpStoc} & 55  & 111 & 166 & 202 & 218 \\  
  & model \eqref{ExpConv} & 0   & 0   & 0   & 0   & 0   \\  
  \hline
  & model \eqref{ExpOnlyDA} & 5.96  & 10.87  & 15.78  & 20.69 & 28.94 \\   
  Inv. cost & model \eqref{ExpStoc} & 7.49  & 13.87  & 20.23  & 27.15 & 35.99 \\  
  & model \eqref{ExpConv} & 7.79  & 14.45  & 21.11  & 28.52 & 38.26 \\  
  \hline
\end{tabular}\label{table:resultsexpansion}\end{table}

As expected, the expansion model~\eqref{ExpOnlyDA} does not
suggest investing in the flexible generating units $\hat{g}_1$ or $\hat{g}_2$, since it assumes that the wind power production is perfectly predictable. Furthermore, notice that, according to
\eqref{ExpConv}, a central planner seeking to meet the renewable
energy target should invest in flexible unit $\hat{g}_1$ (with a
less expensive upward balancing service), but not in flexible unit $\hat{g}_2$ (which offers cheaper downward regulation). The reason for this is that, in order for the system to fully benefit from the cheap downward regulation of $\hat{g}_2$, this unit needs to be dispatched in the day-ahead market out of merit order. This is, however, not possible under the market design considered in expansion model~\eqref{ExpConv}. In contrast, such dispatch decisions are feasible under the market design in expansion model~\eqref{ExpStoc}. Accordingly,
model~\eqref{ExpStoc} suggests investing in unit $\hat{g}_2$ and
not in unit $\hat{g}_1$, with the consequent increase in system
efficiency as downward regulation brings fuel cost savings. Finally, observe that considering forecast errors involves an increase of the installed capacity of wind farm $\hat{w}_1$ in models \eqref{ExpStoc} and \eqref{ExpConv} compared to model \eqref{ExpOnlyDA}, being such capacity higher for the inefficient market design considered in expansion model \eqref{ExpConv}. 

Regarding investment costs, expansion model~\eqref{ExpOnlyDA} unsurprisingly yields the least costly expansion plan, as it supposes an utopian power system operation without wind power forecast errors. The impact of these errors on the investment costs depends on both the underlying market design and the wind energy target. Indeed, for low wind power penetration levels, the efficient market design in expansion model \eqref{ExpStoc} and the inefficient market design in model \eqref{ExpConv} induce expansion plans that are similar in terms of investment costs. However, as the wind energy target is increased, the inefficient market design prompts an expansion plan which is gradually more expensive than the one triggered by the efficient market design. 

\section{Case study}\label{CaseStudy}

In this section we present the results from expansion models 
\eqref{ExpOnlyDA}, \eqref{ExpStoc} and \eqref{ExpConv} in a more realistic 24-bus power system \cite{Grigg99}. The characteristics of the existing generating units are listed in Table \ref{table:24busgendata}. Line susceptances are those in \cite{Grigg99}, while the capacities of all lines are reduced to 175 MW. The demand at each bus is determined according to the parameters provided in \cite{Grigg99} for a peak demand equal to 2850 MW. As in the illustrative example, the demand is assumed to be inelastic, and failures of units and lines are disregarded for simplicity. Forecast errors of demand are likewise neglected.

\begin{table}[htb] \centering \renewcommand*{\arraystretch}{1.5} \setlength{\tabcolsep}{2pt}
\caption{Existing generating units (case study)}
\begin{tabular}{c c c c c c c c c}
  \hline
  {\vphantom{\Large{$A_p$}}} & $x_i/c_{i}$ & $x^{\Plus}_i/c^{\Plus}_i$ & $x^{\Minus}_i/c^{\Minus}_i$ & {\vphantom{\Large{$A_p$}}} & $x_i/c_{i}$ & $x^{\Plus}_i/c^{\Plus}_i$ & $x^{\Minus}_i/c^{\Minus}_i$ \\
  \hline
  $g_1$($n_1$) & 400/25.9 & - &- &  $g_6$($n_{16}$) & 400/19.2  & -   & -  \\
  $g_2$($n_2$) & 575/22.3 & - &- &  $g_7$($n_{18}$) & 120/30.1  & 30/31.1  & 30/29.1  \\
  $g_3$($n_7$) & 500/26.6 & - &- &  $g_8$($n_{21}$) & 100/30.6   & 30/31.6   &30/29.6  \\
  $g_4$($n_{13}$) & 520/21.2   & -   &- &  $g_9$($n_{22}$) & 80/31.1   & 30/32.1   &30/30.1  \\
  $g_5$($n_{15}$) & 475/17.5   & -   &- &  $g_{10}$($n_{23}$) & 450/20.8   & - & -  \\ 
  \hline
\end{tabular}\label{table:24busgendata}\end{table}

In order to increase the penetration of renewable electricity production, the following projects are available:

\begin{itemize}
\item[-] Wind farms of 1000 MW maximum (in blocks of 50 MW) to be located at buses $n_6$, $n_8$, $n_{13}$ and $n_{23}$ with a fixed investment cost $q^0_i=\$25\,000/$year and a variable investment cost of $q_i=\$75\,000/$MW$\cdot$year. 
\item[-] Reinforcement of some existing lines from single to double circuit ($n_6n_{10},n_{11}n_{13},n_{11}n_{14},n_{14}n_{16}$) and the construction of a new line of 51 miles from $n_{12}$ to $n_{21}$ of 175 or 350 MW of capacity and a susceptance of 37.8 p.u.  Fixed and variable investment costs of transmission lines amount to \$16\,400/mile$\cdot$year and \$2.88/MWmile$\cdot$year, respectively. 
\item[-] Six flexible generating units which makes up the installation of one or two additional generating groups to the existing gas-based power plants at buses $n_{18}$, $n_{21}$ and $n_{22}$ with the same characteristics as those provided in Table \ref{table:24busgendata}. Fixed and variable investment costs amount to \$20\,000$/$year and \$25\,000$/$MW$\cdot$year.

\end{itemize}

Wind speed data of 2006 provided by the National Renewable Energy Laboratory (NREL) and corresponding to a site with coordinates 45$^{\circ}$13' N, 96$^{\circ}$55' W is employed in this analysis to model wind power productions at buses $n_6$, $n_8$, $n_{13}$ and $n_{23}$, which are assumed to be perfectly correlated. These data can be freely downloaded from \cite{NREL2010}. A set of 10 scenarios is generated to characterize the variability of demand and stochastic production throughout the planning horizon. Besides, conditional on each scenario at the day-ahead stage, a new set of 10 scenarios modeling the uncertainty of the forecast errors is also generated. Therefore, a total number of 100 scenarios is considered in this study. 

Table \ref{table:resultsexpansioncasestudy} provides expansion planning results for a renewable target of 20\%. Observe that unlike model \eqref{ExpOnlyDA}, expansion models \eqref{ExpStoc} and \eqref{ExpConv} propose the installation of additional flexible generation as well as the line $n_{12}n_{21}$ connecting the area of flexible generation (north) with the area of wind production (east). In addition, it is worth mentioning that the expansion plan suggested by model \eqref{ExpConv} entails an investment cost 7.2\% higher than that of model \eqref{ExpStoc}, which highlights the benefits of an efficient market design in reducing the expansion efforts required to integrate a given amount of renewable production into a power system. 

\begin{table}[htb] \centering \renewcommand*{\arraystretch}{1.5} \setlength{\tabcolsep}{2pt}
\caption{Optimal expansion plans (case study). Capacities in MW and investment cost in m\$.}
\begin{tabular}{l l c c c}
  \hline   
  & & \: Model \eqref{ExpOnlyDA} \: & \: Model \eqref{ExpStoc} \: & \: Model \eqref{ExpConv} \: \\
  \hline
  \multirow{3}{*}{Wind capacity} \: & $n_6$ & 350 & 550 & 350 \\
  & $n_8$ & 500 & 550 & 500 \\
  & $n_{13}$ & 250 & - & 300 \\
  \hline
  \multirow{3}{*}{Flexible Generation} \: & $n_{18}$ & - & 240 & 240 \\
  & $n_{21}$ & - & 160 & 80 \\
  & $n_{22}$ & - & - & 160 \\
  \hline
  \multirow{4}{*}{Line capacity} \: & $n_{6}n_{10}$ & - & 175 & - \\
  & $n_{11}n_{13}$ & - & 175 & 175 \\
  & $n_{14}n_{16}$ & - & - & 175 \\
  & $n_{12}n_{21}$ & - & 350 & 350 \\
  \hline
  \multicolumn{2}{l}{Investment cost} & 132.1 & 143.5 & 153.9 \\
  \hline
\end{tabular}\label{table:resultsexpansioncasestudy}\end{table}

To complement this analysis, we now discuss the implications of disregarding forecast errors by evaluating the expansion plan that results from model \eqref{ExpOnlyDA}. Assuming a market design that efficiently handles forecast errors, the first two rows of Table \ref{table:resultsMC} provide, for different renewable targets, the total expected cost (including both investment and system operating costs in m\$) and the actual level of wind penetration (in parentheses) that result from the expansion plans proposed by models \eqref{ExpStoc} and \eqref{ExpOnlyDA}, respectively. Observe that, although the expansion plan given by model \eqref{ExpOnlyDA} involves a slight increase in the total cost under this market design, the realized wind penetration level is significantly reduced. The reason is that in view of the lack of investments in flexible assets in model (1), an efficient market would reduce the dispatch of stochastic generation to keep the balancing costs low, thus reducing the operational costs, but also the wind penetration level.

Similarly, the third and fourth rows present analogous results considering an inefficient treatment of forecast errors by the market and compare the expansion plans corresponding to models 
\eqref{ExpConv} and \eqref{ExpOnlyDA}. Note that under this type of market design, the expansion plan of model \eqref{ExpOnlyDA} results in a total cost significantly higher than the optimal one, but achieves a wind share level fairly closed to the target. This is due to the fact that an inefficient market would dispatch a high amount of cheap, but uncertain renewable generation and then resort to uneconomical balancing resources (e.g., involuntary load curtailment) to accommodate energy deviations. As a result, the wind target is approximately reached but at an extremely high cost.

\begin{table}[htb] \centering \renewcommand*{\arraystretch}{1.5}\setlength{\tabcolsep}{4pt}
\caption{Impact of forecast errors on expansion planning}
\begin{tabular}{c c c c c}
  \hline
  && \multicolumn{3}{c}{Renewable target}\\
  Market & Expansion & 10\% & 20\% & 30\% \\
  \hline
  \multirow{2}{*}{Effic.}   & model \eqref{ExpStoc} & 416.7(10)  & 451.8(20)  & 489.2(30)   \\ 
  & model \eqref{ExpOnlyDA} & 417.8(9.3)  & 459.5(16.2)  & 498.1(21.8)  \\
  \hline  
  \multirow{2}{*}{Ineffic.}   & model \eqref{ExpConv} & 423.1(10)  & 459.9(20)  & 501.6(30)  \\
  & model \eqref{ExpOnlyDA} & 444.7(9.8)  & 529.0(19.6)  & 602.7(28.5)   \\
  \hline
\end{tabular}\label{table:resultsMC}\end{table}
  
In summary, the results in Table \ref{table:resultsMC} show that, irrespective of how efficiently forecast errors are handled by the market, disregarding these forecast errors when making expansion planning decisions entails undesired outcomes either in terms of system cost or in the achievement of a pre-established renewable target.

\section{Conclusion}\label{Con}

This paper investigates the impact of forecast errors from uncertain renewable generation on the optimal expansion planning of a power system. To this end, we introduce three expansion models that differ in how these errors are handled by the underlying electricity market. Model~\eqref{ExpOnlyDA}, which is the most widespread in the technical literature, simply ignores that these errors occur, while the other two explicitly account for a balancing mechanism to deal with them. Model~\eqref{ExpStoc}, as opposed to \eqref{ExpConv}, represents a market organization where prediction errors are ``perfectly'' managed by co-optimizing the day-ahead and balancing stages. By comparing these three expansion models, we show that: 
\begin{enumerate}
\item Disregarding forecast errors may lead to a highly suboptimal expansion plan.
\item The way this expansion plan is suboptimal, i.e., either in terms of cost efficiency or in terms of renewable energy penetration, is dependent on how efficient the electricity market is in coping with forecast errors.
\item A market design that efficiently handles forecast errors requires lower expansion efforts to integrate a given amount of renewable production into a power system.
\end{enumerate}

As future research, the proposed expansion planning models can be reformulated as multi-year dynamic problems to incorporate timing decisions for expansion projects and long-term uncertainties such as demand growth. Another aspect that requires further investigation is the proposal of solution methods to reduce the computational burden associated with the use of these expansion models in larger power systems. 

% if have a single appendix:
%\appendix[Proof of the Zonklar Equations]
% or
%\appendix  % for no appendix heading
% do not use \section anymore after \appendix, only \section*
% is possibly needed

% use appendices with more than one appendix
% then use \section to start each appendix
% you must declare a \section before using any
% \subsection or using \label (\appendices by itself
% starts a section numbered zero.)
%

%\appendices
%\section{Proof of the First Zonklar Equation}
%Appendix one text goes here.

% you can choose not to have a title for an appendix
% if you want by leaving the argument blank
%\section{}
%Appendix two text goes here.

% use section* for acknowledgement
%\section*{Acknowledgment}

%The authors would like to thank...

% Can use something like this to put references on a page
% by themselves when using endfloat and the captionsoff option.
\ifCLASSOPTIONcaptionsoff
  \newpage
\fi

% trigger a \newpage just before the given reference
% number - used to balance the columns on the last page
% adjust value as needed - may need to be readjusted if
% the document is modified later
%\IEEEtriggeratref{8}
% The "triggered" command can be changed if desired:
%\IEEEtriggercmd{\enlargethispage{-5in}}

% references section
\bibliography{ExpansionPlanningReferences}

% Generated by IEEEtran.bst, version: 1.13 (2008/09/30)
\begin{thebibliography}{10}
\providecommand{\url}[1]{#1}
\csname url@samestyle\endcsname
\providecommand{\newblock}{\relax}
\providecommand{\bibinfo}[2]{#2}
\providecommand{\BIBentrySTDinterwordspacing}{\spaceskip=0pt\relax}
\providecommand{\BIBentryALTinterwordstretchfactor}{4}
\providecommand{\BIBentryALTinterwordspacing}{\spaceskip=\fontdimen2\font plus
\BIBentryALTinterwordstretchfactor\fontdimen3\font minus
  \fontdimen4\font\relax}
\providecommand{\BIBforeignlanguage}[2]{{%
\expandafter\ifx\csname l@#1\endcsname\relax
\typeout{** WARNING: IEEEtran.bst: No hyphenation pattern has been}%
\typeout{** loaded for the language `#1'. Using the pattern for}%
\typeout{** the default language instead.}%
\else
\language=\csname l@#1\endcsname
\fi
#2}}
\providecommand{\BIBdecl}{\relax}
\BIBdecl

\bibitem{EC2009}
E.~E. Commission \emph{et~al.}, ``{Directive 2009/28/EC of the European
  Parliament and of the Council of 23 April 2009 on the promotion of the use of
  energy from renewable sources and amending and subsequently repealing
  Directives 2001/77/EC and 2003/30},'' \emph{Official Journal of the European
  Union Belgium}, 2009.

\bibitem{Bludszuweit2008}
H.~Bludszuweit, J.~Dominguez-Navarro, and A.~Llombart, ``{Statistical Analysis
  of Wind Power Forecast Error},'' \emph{IEEE Trans. Power Syst.}, vol.~23,
  no.~3, pp. 983--991, Aug 2008.

\bibitem{Latorre2003}
G.~Latorre, R.~Cruz, J.~Areiza, and A.~Villegas, ``{Classification of
  publications and models on transmission expansion planning},'' \emph{IEEE
  Trans. Power Syst.}, vol.~18, no.~2, pp. 938--946, May 2003.

\bibitem{Gorenstein1993}
B.~Gorenstin, N.~Campodonico, J.~da~Costa, and M.~V.~F. Pereira, ``{Power
  system expansion planning under uncertainty},'' \emph{IEEE Trans. Power
  Syst.}, vol.~8, no.~1, pp. 129--136, Feb 1993.

\bibitem{Torre2008}
S.~de~la Torre, A.~Conejo, and J.~Contreras, ``{Transmission Expansion Planning
  in Electricity Markets},'' \emph{IEEE Trans. Power Syst.}, vol.~23, no.~1,
  pp. 238--248, Feb 2008.

\bibitem{Baringo2012}
L.~Baringo and A.~Conejo, ``{Transmission and Wind Power Investment},''
  \emph{IEEE Trans. Power Syst.}, vol.~27, no.~2, pp. 885--893, May 2012.

\bibitem{Orfanos2013}
G.~Orfanos, P.~Georgilakis, and N.~Hatziargyriou, ``{Transmission Expansion
  Planning of Systems With Increasing Wind Power Integration},'' \emph{IEEE
  Trans. Power Syst.}, vol.~28, no.~2, pp. 1355--1362, May 2013.

\bibitem{Kirschen2004}
D.~S. Kirschen and G.~Strbac, \emph{{Fundamentals of Power System
  Economics}}.\hskip 1em plus 0.5em minus 0.4em\relax John Wiley \& Sons, 2004.

\bibitem{Pritchard2010}
G.~Pritchard, G.~Zakeri, and A.~Philpott, ``{A Single-Settlement, Energy-Only
  Electric Power Market for Unpredictable and Intermittent Participants},''
  \emph{Oper. Res.}, vol.~58, no. 4-part-2, pp. 1210--1219, 2010.

\bibitem{Morales2013}
J.~M. Morales, M.~Zugno, S.~Pineda, and P.~Pinson, ``{Electricity market
  clearing with improved scheduling of stochastic production},'' \emph{Eur. J.
  Oper. Res.}, vol. 235, no.~3, pp. 765--774, 2014.

\bibitem{birge2011}
J.~R. Birge and F.~Louveaux, \emph{{Introduction to stochastic
  programming}}.\hskip 1em plus 0.5em minus 0.4em\relax Springer, 2011.

\bibitem{dempe2002}
S.~Dempe, \emph{{Foundations of bilevel programming}}.\hskip 1em plus 0.5em
  minus 0.4em\relax Springer, 2002.

\bibitem{Smeers2011}
A.~Ehrenmann and Y.~Smeers, ``\BIBforeignlanguage{English}{{Stochastic
  Equilibrium Models for Generation Capacity Expansion}},'' in
  \emph{\BIBforeignlanguage{English}{{Stochastic Optimization Methods in
  Finance and Energy}}}, ser. {International Series in Operations Research \&
  Management Science}, M.~Bertocchi, G.~Consigli, and M.~A.~H. Dempster,
  Eds.\hskip 1em plus 0.5em minus 0.4em\relax Springer New York, 2011, vol.
  163, pp. 273--310.

\bibitem{Murphy2005}
F.~H. Murphy and Y.~Smeers, ``{Generation Capacity Expansion in Imperfectly
  Competitive Restructured Electricity Markets},'' \emph{Oper. Res.}, vol.~53,
  no.~4, pp. 646--661, 2005.

\bibitem{Nanduri2009}
V.~Nanduri, T.~Das, and P.~Rocha, ``{Generation Capacity Expansion in Energy
  Markets Using a Two-Level Game-Theoretic Model},'' \emph{IEEE Trans. Power
  Syst.}, vol.~24, no.~3, pp. 1165--1172, Aug 2009.

\bibitem{Wogrin2013}
S.~Wogrin, B.~Hobbs, D.~Ralph, E.~Centeno, and J.~Barqu{\'i}n,
  ``\BIBforeignlanguage{English}{{Open versus closed loop capacity equilibria
  in electricity markets under perfect and oligopolistic competition}},''
  \emph{\BIBforeignlanguage{English}{Math. Prog.}}, vol. 140, no.~2, pp.
  295--322, 2013.

\bibitem{floudas1995}
C.~Floudas, \emph{{Nonlinear and Mixed-Integer Optimization: Fundamentals and
  Applications}}.\hskip 1em plus 0.5em minus 0.4em\relax Oxford University
  Press, USA, 1995.

\bibitem{motto2005}
A.~L. Motto, J.~M. Arroyo, and F.~D. Galiana, ``{A mixed-integer LP procedure
  for the analysis of electric grid security under disruptive threat},''
  \emph{IEEE Trans. Power Syst.}, vol.~20, no.~3, pp. 1357--1365, August 2005.

\bibitem{Pineda2014}
S.~{Pineda} and J.~M. {Morales}, ``{Modeling the Impact of Imbalance Costs on
  Generating Expansion of Stochastic Units},'' \emph{ArXiv e-prints}, Feb.
  2014.

\bibitem{Holttinen2005}
H.~Holttinen, ``{Hourly wind power variations in the Nordic countries},''
  \emph{Wind Energy}, vol.~8, no.~2, pp. 173--195, 2005.

\bibitem{Fabbri2005}
A.~Fabbri, T.~Rom{\'a}n, J.~Abbad, and V.~Quezada, ``{Assessment of the Cost
  Associated With Wind Generation Prediction Errors in a Liberalized
  Electricity Market},'' \emph{IEEE Trans. Power Syst.}, vol.~20, no.~3, pp.
  1440--1446, Aug 2005.

\bibitem{Morales2009}
J.~M. Morales, S.~Pineda, A.~J. Conejo, and M.~Carri{\'o}n, ``{Scenario
  Reduction for Futures Market Trading in Electricity Markets},'' \emph{IEEE
  Trans. Power Syst.}, vol.~24, no.~2, pp. 878--888, May 2009.

\bibitem{Grigg99}
C.~Grigg \emph{et~al.}, ``{The IEEE Reliability Test System-1996. A report
  prepared by the Reliability Test System Task Force of the Application of
  Probability Methods Subcommittee},'' \emph{IEEE Trans. Power Syst.}, vol.~14,
  no.~3, pp. 1010--1020, Aug 1999.

\bibitem{NREL2010}
NREL. (2013) {National {R}enewable {E}nergy {L}aboratory. {E}astern Wind
  Dataset.} [Online]. Available at
  ftp://ftp2.nrel.gov/pub/ewits/TimeSeries/LandBased/2006/, 2010.

\end{thebibliography}
\bibliographystyle{IEEEtran}

% biography section
%
% If you have an EPS/PDF photo (graphicx package needed) extra braces are
% needed around the contents of the optional argument to biography to prevent
% the LaTeX parser from getting confused when it sees the complicated
% \includegraphics command within an optional argument. (You could create
% your own custom macro containing the \includegraphics command to make things
% simpler here.)
%\begin{biography}[{\includegraphics[width=1in,height=1.25in,clip,keepaspectratio]{mshell}}]{Michael Shell}
% or if you just want to reserve a space for a photo:

%\begin{IEEEbiography}{Michael Shell}

%\end{IEEEbiography}

% if you will not have a photo at all:
%\begin{IEEEbiographynophoto}{John Doe}
%Biography text here.
%\end{IEEEbiographynophoto}

% insert where needed to balance the two columns on the last page with
% biographies
%\newpage

%\begin{IEEEbiographynophoto}{Jane Doe}
%Biography text here.
%\end{IEEEbiographynophoto}

% You can push biographies down or up by placing
% a \vfill before or after them. The appropriate
% use of \vfill depends on what kind of text is
% on the last page and whether or not the columns
% are being equalized.

%\vfill

% Can be used to pull up biographies so that the bottom of the last one
% is flush with the other column.
%\enlargethispage{-5in}

% that's all folks
\end{document}